\documentclass[11pt,leqno]{amsart}
\usepackage{latexsym,amssymb}
\usepackage{hyperref}
\usepackage{amsmath,amsthm}
\usepackage{amsfonts}
\usepackage[american]{babel}
\usepackage{mathrsfs}
\usepackage{psfrag}
\usepackage{graphicx}
\usepackage{xcolor}

\newcommand{\R}{\mathbb{R}}
\newcommand{\rn}{\mathbb{R}^N}

\newcommand{\bd}{\partial}
\newcommand{\dive}{\text{div}}

\newcommand{\Deltap}{\Delta_{\p}}

\newcommand{\ee}{\textsf{e}}
\newcommand{\ep}{\varepsilon}

\newcommand{\F}{\mathscr{F}}

\newcommand{\Ha}{\mathscr{H}}

\newcommand{\la}{\lambda}
\newcommand{\La}{\Lambda}

\newcommand{\mink}{Minkowski }
\newcommand{\mw}{\mathsf{b}}

\newcommand{\no}{{\mathbf{\nu}}}

\newcommand{\oo}{\Omega}
\newcommand{\ou}{K}
\newcommand{\oua}{\ou_0}
\newcommand{\oub}{\ou_1}
\newcommand{\oula}{\ou_{\la}}
\newcommand{\oum}{\widetilde{\ou}}
\newcommand{\oza}{\oo_0}
\newcommand{\ozb}{\oo_1}
\newcommand{\ozla}{\oo_{\la}}
\newcommand{\ozta}{\oo_{\tau_0}}
\newcommand{\oztb}{\oo_{\tau_1}}
\newcommand{\oztla}{\oo_{\tau_\la}}

\newcommand{\p}{\text{p}}

\newcommand{\s}{\textsf{s}}
\newcommand{\sla}{S_{\la}}

\newcommand{\ttt}{\bar{t}}
\newcommand{\tla}{\tau_{\la}}

\newcommand{\ta}{\tau_0}
\newcommand{\tb}{\tau_1}

\newcommand{\ula}{u_{\la}}
\newcommand{\um}{\widetilde{u}}

\newcommand{\utla}{u_{\tla}}

\newcommand{\uta}{u_0}
\newcommand{\utb}{u_1}
\newcommand{\utta}{u_{\tau_0}}
\newcommand{\uttb}{u_{\tau_1}}

\newcommand{\wn}{\omega_N}

\newcommand{\xx}{\bar{x}}

\numberwithin{equation}{section}
\newtheorem{theorem}{Theorem}[section]
\newtheorem{prop}[theorem]{Proposition}
\newtheorem{cor}[theorem]{Corollary}
\newtheorem{lemma}[theorem]{Lemma}
\newtheorem{remark}[theorem]{Remark}

\newtheorem{definition}[theorem]{Definition}

\begin{document}
\title[Bernoulli problem]{Concavity properties for elliptic free boundary problems}

\author[C. Bianchini]{Chiara Bianchini}
\author[P. Salani]{Paolo Salani}

\address{C. Bianchini, Dip.to di Matematica ``U. Dini'', Universit\`a degli Studi di Firenze, Viale Morgagni 67/A, 50134 Firenze - Italy}
\email{cbianchini@math.unifi.it}

\address{P. Salani, Dip.to di Matematica ``U. Dini'', Universit\`a degli Studi di Firenze, Viale Morgagni 67/A, 50134 Firenze - Italy}
\email{salani@math.unifi.it}
\date{}


\begin{abstract}
We prove some concavity properties connected to nonlinear Bernoulli type free boundary problems. 
In particular, we prove a Brunn-Minkowski inequality and an Urysohn's type inequality for the Bernoulli Constant and we study the behaviour of the free boundary with respect to the given boundary data. 
Moreover we prove an uniqueness result regarding the interior problem.
\end{abstract}

\maketitle

\section{Introduction}
Free-boundary problems of Bernoulli type arise in various physical situations like fluid dynamics, electrochemical machining, optimal insulation and many more.
There are two different kinds of problems, the exterior and the interior one, both concerning nested domains that,
in the classical situation ($\p=2$), represent an annular condenser with a prescribed boundary component while the other one, the free boundary, is  determined (together with a potential function) so that the intensity of the electrostatic field is constant on it.

Throughout the paper, $\Omega$ and $K$, possibly with subscripts, will be domains (bounded connected open sets) in $\rn$, $N\geq 2$, such that $\overline K\subset\Omega$; in fact, they will often be bounded open convex sets. Moreover, if $u\in C^2(\Omega\setminus\overline K)$, we denote by $Du$ and $D^2u$ its gradient and its Hessian matrix, respectively, while, for $\p>1$, we denote by $\Deltap u$ the $\p$-Laplacian of $u$, that is
$$
\Deltap u=\dive(|Du|^{\p-2}Du)\,.
$$

\emph{The exterior Bernoulli problem}

\noindent Given a domain $\ou$ in $\rn$, a real number $\p>1$ and a positive constant $\tau$, the problem consists in looking
for a function $u$ and for a domain $\oo$, containing $\overline\ou$, such that
\begin{equation}\label{Bext}
\begin{cases}
\Deltap u(x)=0 \quad &\text{in }\oo\setminus\overline\ou,\\
u=1 \quad &\text{on }\bd\ou,\\
u=0,\ |Du|=\tau \quad &\text{on }\bd\oo,\\
0<u<1\quad&\text{in }\oo\setminus\overline\ou\, .
\end{cases}
\end{equation}

\bigskip

\emph{The interior Bernoulli problem}

\noindent Given a domain in $\oo\subseteq\rn$, a real number $\p>1$ and a positive constant $\tau$, the problem consists in finding a
function $u$ and a domain $\ou$, contained in $\oo$, such that
\begin{equation}\label{Bint}
\begin{cases}
\Deltap u(x)=0 \quad &\text{in }\oo\setminus\overline\ou,\\
u=0 \quad &\text{on }\bd\oo,\\
u=1,\ |Du|=\tau \quad &\text{on }\bd\ou,\\
0<u<1\quad&\text{in }\oo\setminus\overline\ou\, .
\end{cases}
\end{equation}

\bigskip

If $u$ is a solution of (\ref{Bext}) or (\ref{Bint}), for convenience we will tacitly continue $u$ by $1$ in $K$, throughout the paper.

The boundary condition $|Du|=\tau$ has to be understood in a classical way, in both cases. Regarding the  $\p$-Laplace equations, here we will always consider classical solutions (justified by the convexity of $\Omega$ and $K$, see later).

Notice that, given $\oo$ in (\ref{Bext}) or $\ou$ in (\ref{Bint}), and neglecting the Neumann condition, the function $u$ remains uniquely determined and viceversa.
Hence, we will speak of {\em a solution} of (\ref{Bext}) or (\ref{Bint}) referring indifferently to the sets $\Omega$ and $K$, respectively, or to the corresponding potential function $u$ (or to both); it will be always clear from the context if we are referring to the involved  set or to the corresponding function (or to the couple function-set).

\bigskip

The classical Bernoulli problems regard the case $\p=2$, that is the Laplace operator, and they have been largely investigated since the pioneering work of Beurling \cite{B}. Other references are for instance \cite{A2}, \cite{AC}, \cite{FR} and \cite{F}; see also \cite{CT} and \cite{CM} and references therein. The treatment of the nonlinear case is more recent and mainly due to Henrot and Shahgholian, see for instance \cite{HS1}-\cite{HS4}; see also \cite{AM}, \cite{GK}, \cite{MPS} and references therein.

According to the literature above (see in particular \cite{HS3}),
it is by now well known that, if $\ou$ is convex, a unique classical solution
of the exterior problem exists, for every $\tau>0$; moreover, the
convexity transfers to $\oo$ (and to every level set of $u$) and $\Omega$ is of class $C^{2,\alpha}$.

The interior problem, instead, need not to have a solution for
every domain $\oo$ and for every positive constant $\tau$. However
(see in particular \cite{HS2}), when $\oo$ is convex (with $C^1$ boundary), there exists
a positive constant $\La(\oo)$, named {\em Bernoulli
constant}, such that problem (\ref{Bint}) has no solution if
$\tau<\La(\oo)$, while it has at least one classical solution if
$\tau\ge\La(\oo)$ (and $\ou$ is a $C^{2,\alpha}$ convex set). In \cite{CT} it is proved that, when $\p=2$,
this solution is unique for $\tau=\La(\oo)$; the same property was
not proved for $\p\neq 2$, as far as we know.

In this paper, we consider the convex case and we want to investigate the behaviour of
solutions of the exterior and interior problems with
respect to the data $\tau$ and $\ou$ or $\oo$, respectively.
Moreover, we will prove the uniqueness of the solution
of (\ref{Bint}) corresponding to $\tau=\La(\oo)$ for $\p>1$ and
we will deal with the behaviour of the Bernoulli constant
$\La(\oo)$ and of $(u,\ou)$ with respect to $\oo$.

Our main results are now described in more detail.

\section{Main Results}

Regarding the interior problem, our main result is the following.
\begin{theorem}\label{bmineq}
Let $\oza$ and $\ozb$ be bounded open convex subsets of $\rn$, $N\geq 2$, of class $C^1$.
Let $\lambda\in(0,1)$ and
$$
\ozla=(1-\la)\oza+\la\ozb\,.
$$
Then
\begin{equation}\label{bm}
\La(\ozla)\leq \left[\frac{1-\la}{\La(\oza)} +
\frac{\la}{\La(\ozb)}\right]^{-1}.
\end{equation}
Moreover, equality holds if and only if $\oza$ and $\ozb$ are homothetic.
\end{theorem}

Formula (\ref{bm}) represents a Brunn-Minkowski type inequality for $\La$ and it can be rephrased in the following way: the Bernoulli constant operator $\Lambda:\mathcal{K}\to\R$ is $-1$-convex (i.e. $\Lambda^{-1}$ is concave) in the class $\mathcal{K}$ of bounded convex sets with respect to Minkowski addition (see \S3.1 for definitions).
Notice that, the exponent $-1$ corresponds to the degree of homogeneity of $\La$, that is
\begin{equation}\label{homo}
\La(\alpha\Omega)=\alpha^{-1}\La(\Omega)\quad\text{for every }\alpha>0\,.
\end{equation}
Indeed, $(u(x),\ou)$ is a solution of (\ref{Bint}) in $\Omega$ with boundary condition $|Du|=\tau$ if and only if $(u(\frac x{\alpha}), \alpha\ou)$ solves (\ref{Bint}) in $\alpha\Omega$ with boundary condition $|Du|=\frac 1{\alpha} \tau$.

The proof of Theorem \ref{bmineq}
makes use of a notion of subsolution for problem (\ref{Bint}),
introduced by Beurling \cite{B}, further developed by Acker
\cite{A2} in the case $\p=2$ and then generalized by Henrot and
Shahgholian \cite{HS2} to the case $\p\neq 2$, in combination with
some recent results about the Minkowski addition of quasi-concave
functions (see \cite{CS1} and \cite{LS}).

An almost straightforward and interesting consequence of Theorem \ref{bmineq} is the following property of isoperimetric nature for $\La$ (more appropriately, we should say an Urysohn's type inequality): 
{\em in the class of convex sets with prescribed mean width, the Bernoulli constant attains the minimum value on balls} (for the  definition of mean width of a convex set, see Section \ref{intresult}). 
In other words, we prove the following.
\begin{cor}\label{isop}
Let $\oo$ be a $C^1$ convex domain in $\rn$, $N\geq2$, with mean width $\mw(\oo)=b$.
Let $B$ be a ball with radius equal to $b/2$.
Then
\begin{equation}\label{isopla}
\La(\oo)\geq\La(B)\,,
\end{equation}
and equality holds if and only if $\Omega$ is a ball.
\end{cor}
This result gives an alternative answer to a question posed by Flucher and Rumpf in \cite{FR}: {\em if $\Omega$ is a bounded open convex set, is $\La(\oo)\ge\La(B)$ where $B$ is a ball with the same volume as $\oo$?}

Notice that, due to the Urysohn inequality (\ref{isobV}) and to the monotonicity (\ref{Lamonot}) of $\La$, Corollary \ref{isop} does not imply a positive answer to the Flucher and Rumpf's question, while a positive answer to the latter would imply our result.
Therefore, as far as we now, the question posed in \cite{FR} remains open.
\medskip

Regarding the exterior problem, we prove the following theorem,
that is a concavity type property for the solution of (\ref{Bext}) with respect to $K$ and $\tau$.
\begin{theorem}\label{bernoulliextAB}

Let $\oua, \oub$ be two convex domains in $\rn$ and
$\ta,\tb>0$. \\
Fix $\la\in[0,1]$ and set
$$
\oula=(1-\la)\oua + \la \oub\quad\text{and}\quad\tla=\frac 1{\dfrac{1-\la}{\ta} + \dfrac{\la}{\tb}}\,.
$$
Denote by $(\utta,\ozta)$, $ (\uttb,\oztb)$ and $(\utla,\oztla)$ the solutions of (\ref{Bint}) with data
$(\oua,\ta)$, $(\oub,\tb)$ and $(\oula,\tla)$, respectively, i.e.
\begin{equation}\label{Bext_ab}
\left\{\begin{array}{lll}
\Deltap u(x)=0 \quad &\text{in } \oo_{\tau_i}\setminus \overline{\ou_i}\quad&\\
u=1 \quad &\text{on }\bd\ou_i\quad&\qquad i=0,1,\la\\
u=0,\ |Du|=\tau_i \quad &\text{on }\bd\oo_{\tau_i}\, .&
\end{array}\right.
\end{equation}

Then
\begin{equation}
(1-\la)\ozta + \la\oztb \subseteq \oztla\,,
\end{equation}
and equality holds if and only if $\oua$ and $\oub$ are
homothetic.
\end{theorem}

Theorem \ref{bernoulliextAB} has a counterpart for the interior case that is Proposition \ref{int_concavity} (see  Section 8); this is an extension of \cite[Theorem 1]{CT} and indeed the latter corresponds to the case $\oza=\ozb$, $\p=2$ in our theorem. 
One of the main consequence of the result of Cardaliaguet and Tahraoui is the uniqueness of the solution to the Bernoulli interior problem for $\tau=\La(\oo)$, $\p=2$. 
Following their argument in \cite{CT}, we can extend their result to $\p\neq 2$.
\begin{theorem}\label{unicita}
Let $\oo$ be a $C^1$ convex domain in $\rn$. Then there exists a unique solution to the interior problem (\ref{Bint}) for $\la=\La(\oo)$.
\end{theorem}

The paper is organized as follows. In Section \ref{prelim} we introduce notation and we recall some basic notions. 
In Section \ref{lem} we prove a monotonicity property for the norm of the gradient of solutions to $\p$-Laplace equation in convex rings, a technical result which has its own interest.
In Section \ref{extresult} we prove Theorem \ref{bernoulliextAB}. Section \ref{intresult} is devoted to the proof of
Theorem \ref{bmineq} and Corollary \ref{isop}. Finally, in Section \ref{Sectunicita} we prove Proposition \ref{int_concavity} and Theorem \ref{unicita}.

\section{Notation and Preliminaries}\label{prelim}

In the $N$-dimensional Euclidean space, $N\geq 2$, we denote by $\langle\cdot,\cdot\rangle$ the classical Euclidean scalar product and by $|\cdot|$ the Euclidean norm. 
For $K\subseteq\rn$, we denote by $\overline K$ its closure and by $\bd K$ its boundary. 
$\Ha^m$ indicates the $m$-dimensional Hausdorff measure. By $B$ we denote the unit ball in $\rn$, that is $B=\{x\in\rn\,:\,|x|< 1\}$. 
Moreover we set $\omega_N=\Ha^N(B)$ and
$$
S^{N-1}=\bd B=\{x\in\rn\,:\,|x|= 1\};
$$ 
hence $\Ha^{N-1}(S^{N-1})=N\omega_N$.

    \subsection{\mink addition and support function of convex sets}
Let $K$ be a subset of $\rn$ and let $\alpha\ge 0$; we set $\alpha K = \{ \alpha x :\ x\in K \}.$
Let $K_0,K_1\subseteq\rn$; we define their \mink sum $K_0+K_1$ as
$$
K_0+K_1 = \{ x_0+x_1\,:\, x_0\in K_0,\, x_1\in K_1 \}.
$$
Let $\lambda=(\lambda_1,\dots,\lambda_m)$, such that $\sum_{i=1}^m\la_i=1$, $\la_i\ge 0$, for every $i=1,...,m$, and let $K_1,\dots,K_m\subseteq\rn$; we set
$$
K_\lambda=\sum_{i=1}^m\lambda_iK_i=\left\{\sum_{i=1}^m\lambda_ix_i\,:\,x_i\in K_i,\, i=1,\dots,m\right\}\,.
$$
Notice that, if $K_1,\dots,K_m$ are convex sets, then $K_\lambda$ is convex as well.

The support function $h_K(\cdot)=h(K,\cdot) : \rn\to [0,+\infty)$, of a convex set $K$ is defined in the following way:
$$
h_K(v)= \sup_{x\in K} \langle v,x \rangle, \qquad v\in\rn.
$$
By definition $h_K$ is obviously homogeneous of degree one and, as supremum of linear function, it is convex.
Moreover, for every $a\geq 0$ and every $K,L$ convex sets, it holds
\begin{equation}\label{nuova}
\begin{array}{ll}
& h_{aK} = ah_K,\\
& h_{K+L} = h_K + h_L.
\end{array}
\end{equation}
We refer to \cite{S} for more details and properties of convex sets and support functions.

    \subsection{Quasi-concave and $Q^2_-$ functions}
An upper semicontinuous function $u:\rn\to\R\cup\{\pm\infty\}$ is said quasi-concave if it has convex superlevel sets, or,
equivalently, if
$$
u\left( (1-\la)x_0+\la x_1 \right)\ge \min\{ u(x_0),u(x_1) \},
$$
for every $\la\in[0,1]$, and every $x_0, x_1\in\rn$.
If $u$ is defined only in a proper subset $\Omega$ of $\R^n$, we extend $u$ as $-\infty$ in $\R^n\setminus\Omega$ and we say that $u$ is
quasi-concave in $\Omega$ if such an extension is quasi-concave in $\rn$.
In an analogous way, $u$ is quasi-convex if $-u$ is quasi-concave, i.e. if it has convex sublevel sets.
Obviously, if $u$ is concave (convex) then it is quasi-concave (quasi-convex).

A special subclass of quasi-concave functions was introduced and studied in \cite{LS}.
\begin{definition}
Let $u$ be a function defined in an open set $\Omega\subset\R^n$; we say that  $u$ is a $Q^2_-$ function at a point $x\in \Omega$ (and we write $u\in Q^2_-(x)$) if:
\begin{enumerate}
\item $u$ is of class $C^2$ in a neighborhood of $x$;
\item its gradient does not vanish at $x$ (i.e. $|Du(x)|>0$);
\item the principal curvatures of $\{ y\in\R^n\ |\ u(y)= u(x)\}$ with respect to the normal $-\frac{Du(x)}{|Du(x)|}$ are positive at $x$.
\end{enumerate}
\end{definition}
In other words, a $C^2$ function $u$ is $Q^2_-$ at a regular point $\xx$ if its level set $\{x\,:\,u(x)=u(\xx)\}$ is a regular convex surface (oriented according to $-Du$), whose Gauss curvature does not vanish in a neighborhood of $\xx$.
By $u\in Q^2_-(\Omega)$ we mean $u\in Q^2_-(x)$ for every $x\in\Omega$.

    \subsection{The support function of a $Q^2_-$ function}
Since a continuous function is completely known if one knows all its level sets, and since every compact convex set is univocally  determined by its support function, we can  associate to every quasi-concave function $u$, a function $h:\rn\times\R\to\R\cup\{\pm\infty\}$, such that, for every fixed $t\in\R$, $h(X,t)$ is the support function of the superlevel set $\{x\in\rn\,:\,u(x)\ge t\}$ evaluated at $X$.
We will refer to $h$ as the support function of the function $u$. Notice that $h$ obviously depends on $u$; sometimes we will stress such a dependence by writing $h_u$, but in general we will avoid this and we will use the subscript with $h$ to indicate partial differentiation.

Next we recall some properties of $Q^2_-$ functions and their support functions from \cite{LS}.
Let $\Omega$ be a convex domain and $u$ be a $Q^2_-(\Omega)\cap C(\overline\Omega)$ function such that $u=0$ on $\partial\Omega$ and $\max_{x\in\overline\Omega}u(x)=1$.
For every $t\in(0,1)$ and for every $X\in\R^N\setminus\{0\}$, there exists a unique point $x(X,t)$ such that
\begin{equation}\label{defxt}
u(x(X,t))=t,           \qquad      -\,\frac{Du(x(X,t))}{|Du(x(X,t))|} = \frac{X}{|X|}.
\end{equation}
In fact $x(X,t)$ is the unique point on $U(t) =\{x\in\rn\ :\ u(x) = t \}$ such that
$$
\langle x(t);X \rangle =\max_{y\in U(t)}\langle y;X \rangle.
$$

Moreover, due to the $C^1$ regularity of $u$, $h\in C^2((\rn\setminus\{0\})\times (0,1))$ with
\begin{equation}\label{ht<}
h_t(X,t)<0\qquad\text{ for every }X\in\rn\setminus\{0\},\,t\in(0,1),
\end{equation}
and
\begin{equation}\label{LS2}
Du(x(X,t)) = -\frac {X}{h_t(X,t)}\,.
\end{equation}
(Here and later $h_t$ denotes the derivative of $h$ with respect to the $t$ variable)

The following expression of the $\p$-Laplacian of a $Q^2_-$ function in terms of $h$ has been proved (for $|X|=1$) in \cite[Proposition 1]{CS1} and it can be also easily deduced from \cite{LS}:
\begin{equation}\label{LS3}
\Deltap u = \dfrac 1{(-h_t)^{\p-1}} \left( (\p-1)h_{tt} - h_t^2 C
-(\p-1)\sum_{i=1}^{N-1}k_ih_{ti}^2 \right)\,,
\end{equation}
where $C$ denotes the mean curvature of the level set $\{u=t\}$ and the left-hand side is evaluated at $x(X,t)$, while the right-hand side is calculated at $(X, t)$ or, equivalently, the left-hand side is evaluated at $x$ while the right-hand side is calculated at $(\frac{Du(x)}{|Du(x)|},u(x))$; here $X$ must be a unitary vector, i.e. $|X|=1$.

    \subsection{\mink addition of functions}
Let $u_1,\dots,u_m$ be upper semicontinuous functions defined in $\Omega_1\dots\Omega_m\subset\rn$, respectively, and let $\la=(\la_1,\dots,\la_m)$ be such that $\sum_{i=1}^m\la_i = 1$, $\la_i\ge 0$, $i=1,...,m$.
The \mink linear combination of the functions $u_i$ with ratio $\la$ is the upper semicontinuous function $\ula$ whose super-level sets $U_{\la}(t)=\{u_\lambda\geq t\}$ are the \mink linear combination of the corresponding super-level sets $U_i(t)=\{u_i\geq t\}$ of $u_i$, i.e.
$$
U_\lambda(t)=\sum_{i=1}^m\la_iU_i(t), \quad\text{ for every }\quad t\in\R,
$$
and
$$
u_\lambda(x)=\sup\{t\,:\, x\in U_\lambda(t)\}\,.
$$
Notice that this operation preserves the quasi-concavity; in particular the class of $Q^2_-$ functions is closed with respect to Minkowski addition,
that is:
{\em let $u_i\in Q^2_-(\Omega_i)$ for $i=1,\dots,m$, then $\ula \in Q^2_-(\Omega_{\la})$} (see \cite{CS1}, \cite{LS}).

Moreover, by (\ref{nuova}), it holds
$$
h_{(1-\la)u_0+\la u_1} = (1-\la) h_{u_0} + \la h_{u_1}.
$$

\section{An auxiliary lemma}\label{lem}
\begin{lemma}\label{DucrC2+}
Let $\oo\setminus\overline\ou$ be a bounded convex ring (i.e. $\oo$ and $\ou$ are bounded convex domains with $\overline{\ou}\subseteq \oo$), and let $u\in C^2(\oo\setminus\overline\ou)\cap C(\overline{\oo\setminus\ou})$ solve
\begin{equation}\label{eqne}
\begin{cases}
\Deltap u(x) \ge 0 \quad &\text{in }\oo\setminus\overline\ou,\\
u=1 \quad &\text{on }\bd\ou,\\
u=0 \quad &\text{on }\bd\oo\, .
\end{cases}
\end{equation}
Moreover, assume that $u\in Q^2_-(\oo)$.
Then for every fixed direction $\theta\in S^{N-1}$, $|Du(x(\theta,t))|$ is strictly increasing with respect to $t\in(0,1)$.
\end{lemma}
\begin{proof}
For simplicity reasons we denote by $x(t)$ the point $x(\theta,t)$ defined by (\ref{defxt}).
Notice that $x(t)$ is a regular curve for $t\in[0,1]$, since $u$ is a $Q^2_-$ function.
By definition of $x(t)$, $\theta$ is the outer unit normal vector to the level set $U(t)=\{x\in\rn\ :\ u(x)=t \}$ at $x(t)$; hence $Du(x(t))$ is parallel to $\theta$ and it points in the opposit direction.

By assumption, $|Du(x(t))|>0$ and all the principal curvatures of $U(t)$ at $x(t)$ are positive for every $t\in(0,1)$ and hence the mean curvature $C$ of $U(t)$ is positive at $x(t)$.

Using the fact that $\Deltap u \ge0$ in $\oo\setminus\overline\ou$, together with (\ref{ht<}) and the positivity of $C$, formula (\ref{LS3}) entails
$$
h_{tt}\ge \frac{C}{\p-1}h_t^2 + \sum_{i=1}^{N-1}k_ih_{ti}^2 > 0.
$$
Hence by (\ref{ht<})-(\ref{LS2}) we obtain
$$
\frac{d}{dt}|Du(x(t))| = \frac{d}{dt} \frac{-1}{h_t(\theta,x(t))} = \frac{h_{tt}(\theta,x(t))}{h_t^2(\theta, x(t))} > 0,
$$
which implies the stated monotonicity of the gradient, with respect to the level parameter.
\end{proof}

\begin{remark}
Notice that in the above proof we just used $C>0$. Hence the $Q^2_-$ assumption can be weakened in the following way:
$|Du|\neq 0$ and every level set has always positive mean curvature.
\end{remark}

The previous Lemma is an infinitesimal version of \cite[ Lemma 2.2 ]{HS3} and it has its own interest; also compare it with 
\cite[Lemma 2.7]{HS4}, which contains a similar result for the flow curves of $Du$.

\section{Proof of Theorem \ref{bernoulliextAB}}\label{extresult}

Let $\sla = (1-\la)\ozta + \la\oztb$.
In order to show that $\sla \subseteq \oztla$, we compare the {\mink} linear combination $\ula$ of the functions $\utta, \uttb$, with $\utla$.

Notice that $\utta\in Q^2_-(\ozta)$ and $\uttb\in Q^2_-(\oztb)$ by \cite[Theorem 1]{L}, since  $\Omega_{\tau_i}$, $i=0,1$, must be convex (see \cite{HS1}, \cite{HS3}, \cite{GK}, for instance). Then $\ula\in Q^2_-(\ozla)$ (see \cite{CS1}, \cite{LS}).

Let us indicate by $V(t)$ the superlevel sets of $\ula$ of level $t\in[0,1]$, that is
$$V(t)=\{ x\in\rn\,:\,\ula(x)\ge t \}=(1-\la)\{x\in\rn\,:\,\utta(x)\geq t\}+\la\,\{x\in\rn\,:\,\uttb(x)\geq t\}\,.$$
Notice that $V(t)$ is convex for every value of $t$.

Set
$$
\ttt = \inf \{ t\in[0,1]\ :\ \overline{V(t)}\subseteq\overline{\oztla} \}\,,
$$
and, by contradiction, assume $\ttt>0$, since the case $\ttt=0$ easily implies the thesis.

By the regularity of $\ula$, the infimum in the definition of $\ttt$ is in fact a minimum and there exists at least one point $\xx\in \bd V(\ttt)\cap\bd\oztla$ (while $V(\ttt)\subset\overline\oztla$).
Then the outer unit normal vectors to $V(\ttt)$ and to $\oztla$ at $\xx$ coincide; let us denote this vector by $\no$. Since
$\bd V(\ttt)$ and $\bd\oztla$ are level sets of $\ula$ and $\utla$, respectively, we have
\begin{equation}\label{Dunu}
\dfrac{D\ula(\xx)}{|D\ula(\xx)|}=\dfrac{D\utla(\xx)}{|D\utla(\xx)|}=-\no\,.
\end{equation}
By definition of \mink sum, there exist $x_0\in\bd\ozta(\ttt)$, $x_1\in\bd\oztb(\ttt)$ such that $\xx=(1-\la)x_0+\la x_1$ with
$$
\dfrac{D\utta(x_0)}{|D\utta(x_0)|} = \dfrac{D\uttb(x_1)}{|D\uttb(x_1)|} = -\no\,;
$$
then (see \cite{CS1},\cite{LS})
\begin{equation}\label{LS}
|D\ula(\xx)| = \left( \dfrac{1-\la}{|D\utta(x_0)|} + \dfrac{\la}{|D\uttb(x_1)|} \right)^{-1}.
\end{equation}
Notice that, since $\ttt>0$, Lemma \ref{DucrC2+} yields
\begin{equation*}
|D\utta(x_0)|>\tau_0\,,\quad |D\uttb(x_1)|>\tau_1\,.
\end{equation*}
Hence
(\ref{LS}) gives
\begin{equation}\label{Dumaggioretau}
|D\ula(\xx)| > \tla\,.
\end{equation}
Now set
$$
w=\ula-\utla.
$$
By \cite{CS1}, \cite{LS} it holds $\Deltap \ula \ge 0$  in $\sla\setminus\overline\oula,$ with $\ula=1$ on $\bd \oula$, $\ula = 0$ on $\bd \sla$.
Hence 
by the comparison principle $u_\lambda\le u_{\tau_\lambda}+{\overline t}$ so that
$$
w(x)\le \ttt \text{ in }V(\ttt)\setminus \overline\oula;
$$
then
$$
\max_{\overline{V(\ttt)\setminus\oula}}w=\max_{\bd (V(\ttt)\setminus\oula)} w = \max_{\bd V(\ttt)} w = \ttt=w(\xx)\,.
$$
As $\ttt>0$, $w$ is not constant on $V(\ttt)$ and the Maximum Principle gives
$$
\max_{\overline{V(\ttt)\setminus\oula}}w=\max_{\bd (V(\ttt)\setminus\oula)} w = \max_{\bd V(\ttt)} w = \ttt=w(\xx),
$$
since $w$ vanishes on $\bd \oula$ and $\utla\ge 0$ on $\bd V(\ttt)$ with $\utla(\xx) =0$ (for $\utla = 0$ on $\bd\oztla$).

In particular, $\xx$ is an absolute maximum  point, hence
\begin{equation}\label{wle}
\dfrac{\bd w(\xx)}{\bd\no} \ge 0.
\end{equation}
On the other hand, by (\ref{Dunu}) and by the definition of $w$, it holds
$$
Dw(\xx)= D\ula(\xx)-D\utla(\xx)=\left(|D\utla(\xx)|-|D\ula(\xx)|\right)\no\,,
$$
whence, by the definition of $\utla$ and by (\ref{Dumaggioretau}),
$$
\dfrac{\bd w(\xx)}{\bd\no}=\langle Dw(\xx),\no\rangle=|D\utla(\xx)|-|D\ula(\xx)|=\tla-|D\ula(\xx)|<0\,,
$$
which contradicts (\ref{wle}).

This shows that it is not possible to assume $\ttt>0$ and hence $V(0)= \sla\subseteq \oztla$.
\bigskip

Let us now consider the equality case.
If $\oua,\oub$ are homothetic, then it is enough to notice that $(u(x),\oo)$ is the solution to the problem corresponding to some $\ou$ and $\tau$ if and only if $(u(\frac x{\alpha}),\alpha\oo)$ is the solution corresponding to $\alpha\ou$ and $\frac{\tau}{\alpha}$.

On the other hand, assume that $\sla=\oztla$; then the functions $\ula$ and $\utla$ coincides on $\bd\oztla\cup\bd\oula=\bd(\oztla\setminus\overline\oula)$ and $|D\ula|=|D\utla|$ on $\bd\oztla$.
Using the Maximum Principle and the Hopf Lemma (see \cite{Tolk}), $\ula$ and $\utla$ must coincides in $\oztla\setminus\overline\oula$.
Hence $\Deltap\ula=0$ in $\oztla\setminus\overline{\oula}$ which implies that the corresponding level sets of $\utta$ and $\uttb$ are homothetic (see the proof of the equality case in the Brunn-Minkowski inequality for the $\p$-capacity in \cite{CS1}). Hence $\oua$ and $\oub$ are homothetic domains.

\section{The Brunn-Minkowski inequality and an isoperimetric inequality for the Bernoulli constant}\label{proofbm}\label{intresult}

Before going into detail of the proof of Theorem \ref{bmineq}, we need to recall the notion of subsolution and maximal solution of problem (\ref{Bint}).

Let $\Omega$ be a subset of $\rn$; $\F(\Omega,\tau)$ is the class of functions $v$ that are Lipschitz continuous on $\overline{\Omega}$ and such that
\begin{equation*}
\begin{cases}
\Deltap v \ge 0 &\qquad\text{in } \{ v<1\}\cap \Omega\\
v =0           &\qquad\text{on }\bd \Omega\\
|Dv| \le \tau  &\qquad\text{on }\bd\{v<1 \}\cap \Omega\,.
\end{cases}
\end{equation*}
If $v\in\F(\Omega,\tau)$ we call it a {\em subsolution}.

With abuse of terminology and notation, we say that
a set $K$ is a {\em subsolution}, and we possibly write $K\in\F(\Omega,\tau)$
or $(v,K)\in\F(\Omega,\tau)$,
if $K=\{x\in\Omega\,:\,v(x)\geq 1\}$ for some $v\in\F(\Omega,\tau)$.

Essentially, $v$ and $\ou$ are subsolutions if $v$ solves
$$
\begin{cases}
\Deltap v\ge 0       &\qquad\text{ in }\oo\setminus\overline{\ou}\\
v=0                  &\qquad\text{ on }\bd\oo\\
v=1,\ |Dv|\le \tau   &\qquad\text{ on }\bd\ou\,.
\end{cases}
$$
In \cite{HS2} Henrot and Shahgholian proved that, when $\oo$ is convex, the Bernoulli constant can be characterized in the following way
\begin{equation}\label{bernconstant}
\La(\oo)= \inf\{ \tau\,:\,\F(\oo,\tau) \neq\emptyset\}\,,
\end{equation}
and a solution to (\ref{Bint}) exists if and only if $\tau\geq\La(\oo)$. 
In such a case, they proved, in particular, that there exist  a {\em largest set}
\begin{equation}\label{tildak}
\tilde K(\oo,\tau)=\cup_{A\in\F(\oo,\tau)}A\,,
\end{equation}
and a {\em maximal function} $\tilde u=\sup_{v\in\F(\oo,\tau)}v$, such that $(\tilde u,\tilde K)\in\F(\oo,\tau)$; in fact,  the couple $(\tilde u,\tilde K)$ solves (\ref{Bint}), and it is called {\em maximal solution}. Moreover, the set $\tilde K$ is convex and the function $\tilde u$ is quasi-concave.

Notice that, from (\ref{bernconstant}) the monotonicity of $\La$ with respect to the inclusion easily follows:
\begin{equation}\label{Lamonot}
\text{ if }\quad\Omega_A \subseteq \Omega_B \quad\text{ then }\quad \La(\Omega_A) \ge \La(\Omega_B)\,.
\end{equation}

Now we can proceed with the proof.
\begin{proof}[Proof of Theorem \ref{bmineq}]
We denote by $(u_i,\ou_i)$, $i=0,1$, the maximal solutions of (\ref{Bint}) corresponding to $\Omega_i$ with $\tau=\La(\Omega_i)$, $i=0,1$ (see Theorem \ref{unicita} for a proof of the uniqueness of the solution when $\tau=\La(\oo)$).
By \cite{HS2} and \cite{L}, $u_i\in Q^2_-(\oo_i\setminus\overline K_i)$, $i=0,1$.

Set $\sla = (1-\la)\oua + \la \oub\,$ and
$$
\tla = \frac 1{\frac{1-\la}{\La(\oza)} + \frac{\la}{\La(\ozb)}}.
$$
Let $\ula$ be the \mink addition of the functions $\uta, \utb$ with ratio $\la$; then $\ula\in Q^2_-(\ozla\setminus\overline{\sla})$ and it solves (see \cite{CS1}, \cite{LS})
\begin{equation}\label{pbula}
\begin{cases}
\Deltap \ula(x) \geq0     &\qquad \text{ in }\ozla\setminus\overline\sla\\
\ula=0                    &\qquad \text{ on }\bd\ozla\,,\\
\ula=1, |D\ula| = \tla    &\qquad \text{ on }\bd\sla\,.
\end{cases}
\end{equation}

This proves that $(\ula,\sla)\in\F(\ozla,\tla)$ and hence $\La(\ozla)\le \tla$ by (\ref{bernconstant}).
\bigskip

Now, let us turn to the equality case.
If $\oua,\oub$ are homothetic, then equality holds by the homogeneity of $\La$.

On the other hand, if $\La(\ozla) = \tla$, consider again the Minkowski linear combination $\ula$ of $u_0$ and $u_1$; then, by (\ref{pbula}), $\ula\in\F(\ozla,\La(\ozla))$.

Let $v$ be the $\p$-capacitary function of $\ozla\setminus\overline\sla$:
$$
\begin{cases}
\Deltap v(x) =0 &\qquad \text{ in }\ozla\setminus\overline\sla\\
v=0             &\qquad \text{ on }\bd\ozla\\
v=1             &\qquad \text{ on }\bd\sla\,.
\end{cases}
$$
Thanks to the convexity of $\ozla$ and $\sla$, $v$ has convex super-level sets (see \cite{G}, \cite{Ka}, \cite{L}).

If $\ula\neq v$, the Hopf Lemma (see \cite{Tolk}) gives
$$
\tla=|D\ula|>|Dv| \qquad\text{ on }\bd\sla\,,
$$
hence
$$
\tau=\max_{x\in\bd\sla}|Dv|<\tla\,.
$$
This implies that $v\in\F(\ozla,\tau)$ with $\tau<\La(\ozla)$, which contradicts (\ref{bernconstant}). Hence $\ula$ coincides with $v$, which entails that all the corresponding level sets of $\uta,\utb$ are homothetic (see \cite{CS1}).
\end{proof}

An isoperimetric type inequality for $\La$ easily follows from Theorem \ref{bmineq}. Before proving Corollary \ref{isop}, let us discuss a little bit this result.

In \cite{HS2}, Henrot and Shahgholian considered the class of sets with fixed minimum width (i.e. such that the diameter of the largest ball inscribed is fixed). They proved that for every convex set $K$, with minimum width $d=2R$, it holds
\begin{equation}\label{isopminwid}
\La(\oo)\ge \frac 1R.
\end{equation}
Notice that equality in (\ref{isopminwid}) does not hold even when $\Omega$ is a ball, since the Bernoulli constant of a ball $B_R$ of radius $R$ is
computed as
\begin{equation}\label{LaBall}
\La(B_R)=
\begin{cases}
\left( \dfrac{N-1}{\p-1} \right)^{\frac{N-1}{N-\p}} \dfrac 1R  &\qquad\text{ if } N\neq \p;\\
&\\
\dfrac {\ee}R &\qquad\text{ if }N=\p.
\end{cases}
\end{equation}

Another estimate simply follows from the monotonicity (\ref{Lamonot}) of $\La$: for every convex set $K$ with outer radius $d$ (i.e. such that the radius of the smallest ball containing $K$ is fixed equal to $d$), it holds
\begin{equation}\label{outerradius}
\La(K) \ge \La(B_d)\,.
\end{equation}


Here we consider the class of convex domains with fixed mean width.
The mean width $\mw (\Omega) $ of a convex set $\Omega$ is defined as
$$
\mw (\Omega) = \dfrac{2}{N\wn} \int_{S^{N-1}} h_\Omega(\theta)\, d\Ha^{N-1}(\theta).
$$
We recall that the following Urysohn's inequality holds in the class of convex sets:
\begin{equation}\label{isobV}
\frac{V(\Omega)}{\omega_n}\leq\left(\frac{\mw(\Omega)}{2}\right)^n\,.
\end{equation}

Trivially it follows that the mean width of a convex set is less or equal then twice the outer radius (and equality holds only for balls).
Hence inequality (\ref{isopla}) is stronger than (\ref{outerradius}).

\begin{proof}[Proof of Corollary \ref{isop}]
Let $\Omega$ be a subset of $\R^N$ with mean width $b$ and {\em Steiner point} $s$.
We recall that the Steiner point $\s(\Omega)$ of a convex set $\Omega$ can be defined as
$$
\s(\Omega) = \dfrac 1{\wn} \int_{S^{N-1}} \theta\, h_\Omega(\theta)\ d\Ha^{N-1}(\theta).
$$
By Hadwiger's Theorem (see \cite{S}, Section 3.3) there exists a sequence of rotations $\{\rho_n\}$ such that
$$
\Omega_n = \frac 1n (\rho_1\Omega+...+\rho_n\Omega),
$$
converges, in the Hausdorff metric, to a ball $B$.
Notice that, since the mean width is invariant under rigid motions and it is \mink additive (see \cite{S}, Section 1.7), $\mw(\Omega_n)=b$ for every $n$ and hence $\mw(B)=b$. Moreover $\s(\Omega_n)=s$ for every $n$ for the same reasons and hence $B$ is the ball with radius $R=\frac b2$ and center $s$.

By Theorem \ref{bmineq} it holds
$$
\La(\Omega_n) \le \La(\Omega),
$$
since $\La(\rho\Omega) = \La(\Omega)$ for any rotation $\rho$.

Since $\Omega_n$ converges to $B$  in the Hausdorff metric as $n$ tends to infinity, there exists $m$ such that $\Omega_n\subseteq B_{R+\frac 1n}$ for every $n\ge m$, where $B_{R+\frac 1n}$ is the ball with radius $R+\frac 1n$ and center $s$.
Then, by (\ref{Lamonot}) and (\ref{LaBall}),
$$\La(\Omega_n)\ge\La(B_{R+\frac 1n}) = \left(\frac{|N-1|}{|\p-1|}\right)^{\frac{N-1}{N-\p}} \frac 1{R+\frac 1n}\,,$$
which converges to $\La(B)$ as $n$ tends to infinity, and this complete the proof of (\ref{isopla}).

Let us now characterize the minimizers; in particular, let us show that balls are the only ones.
Since the Bernoulli constant is invariant under translations, we can consider convex sets with assigned mean width and also fixed Steiner point; hence there exists a unique ball $B$ in the class.
By contradiction, assume that $\Omega$ belongs to this class with $\La(\Omega)=\La(B)$ and $\Omega$ does not coincide with $B$.
Then by Theorem \ref{bmineq} and (\ref{homo}) we have
$$
\La\left(\frac{\Omega+B}2\right)<\left( \frac 1{2\La(\Omega)} + \frac 1{2\La(B)} \right)^{-1}= \La(B),
$$
with
$$
\mw\left(\frac{\Omega+B}2\right)=\mw(B),
$$
by the \mink additivity of the mean width. This contradicts (\ref{isopla}).
\end{proof}

\begin{remark}
The inequality stated in Theorem \ref{bmineq} is not surprising; indeed analogous inequalities hold for several other functionals from calculus of variations. 
In particular, the proof of Corollary \ref{isop} works almost unchanged every time a Brunn-Minkowski type inequality holds. 
Hence it can be adapted, for instance, to the cases of the first Dirichlet eigenvalue of the $\p$-Laplacian or to the eigenvalue of the Monge-Amp\`ere (see \cite{CCS} and \cite{Sa} for the related Brunn-Minkowski inequalities). 
In the latter case, a direct proof of this inequality via an ad hoc symmetrization procedure is given in \cite{Tso}, where also stronger inequalities for the eigenvalues of the other Hessian equations are also proved. On the other hand, the Monge-Amp\`ere eigenvalue is the only case where the inequality involving mean width is known to be optimal (see again \cite{Tso}). 
In particular, the Faber-Krahn inequality for the first Dirichlet eigenvalue $\lambda_{\p}$ of the $\p$-Laplacian (see \cite{KF} and references therein), asserts that balls give the minimum value of $\lambda_{\p}(\Omega)$ among all sets $\Omega$ with given volume. 
Hence we believe that inequality (\ref{isopla}) is not optimal and we conjecture that a positive answer to the Flucher and Rumpf question can be given.
\end{remark}

\section{Uniqueness for the interior problem}\label{Sectunicita}

 Before proving Theorem \ref{unicita}, we prove a connected concavity property for maximal solutions of (\ref{Bint}), which is a generalization of
\cite[Theorem 1.1]{CT}.
\begin{prop}\label{int_concavity}
Let $\oza$ and $\ozb$ be two $C^1$ convex domains in $\rn$ and let $\ta\ge\La(\oza)$, $\tb\ge\La(\ozb)$.
Set
$\lambda\in(0,1)$, $\ozla=(1-\la)\oza+\la\ozb$ and
$$
\tla = \frac 1{\frac{1-\la}{\ta}+\frac{\la}{\tb}}.
$$
Then
$$
(1-\la)\oum(\oza,\ta) +\la \oum(\ozb,\tb) \subseteq
\oum(\ozla,\tla),
$$
where $\tilde K(\oo_i,\tau_i)$ denotes the largest set for $i=0,1,\la$.
\end{prop}
\begin{proof} Following the argument in the proof of Theorem \ref{bmineq}, let us consider
$$
\sla=(1-\la)\oum(\oza,\ta)+\la\oum(\ozb,\tb),
$$
and $\ula$ the \mink addition of the corresponding maximal functions $\um_0,\um_1$.
Since $\ula$ solves problem (\ref{pbula}), $(\ula,\sla)\in\F(\ozla,\tla)$ which, by (\ref{tildak}), implies
$$
\sla \subseteq \oum(\ozla,\tla).
$$
\end{proof}

\begin{proof}[Proof of Theorem \ref{unicita}]
Our proof essentially follows the argument of \cite[Theorem 1.2]{CT}.

Let us indicate by $(u,K)$ the maximal solution to problem (\ref{Bint}) with $\tau=\La(\oo)$.
Assume by contradiction that there exists another solution $(v,D)$ of problem (\ref{Bint}) with $\tau=\La(\oo)$, and assume $D\neq K$ (hence $D\subsetneq K$).
Consider $D^*$ the convex hull of $D$ and let $v^*$ be the $\p$-capacitary function of $D^*$ with respect to $\oo$, that is the solution to the Dirichlet problem:
$$
\begin{cases}
\Deltap v^*=0\qquad&\text{ in }\oo\setminus D^*\\
v^*=0      \qquad&\text{ on }\bd\oo\\
v^*=1        \qquad&\text{ on }\bd D^*.
\end{cases}
$$
By \cite{HS2}, Lemma 2.5, $|Dv^*|\le\La(\oo)$, that is $v^*\in\F(\oo,\La(\oo))$.

Notice that $D^*\subseteq K$, but it cannot coincide with $K$.
Indeed, by contradiction, assume that $D^*=K$, then there should exist $\xx\in\bd D\cap\bd K$; by Hopf's Lemma, using the fact that $D\subsetneq K$, we obtain $|Dv(\xx)| > |Du(\xx)|$, which contradicts $|Dv|=|Du|=\La(\oo)$ on $\bd K\cap\bd D$.

Now, fix $\la\in(0,1)$, let $\sla= (1-\la)K + \la D^*$ and let $\ula$ be the {\mink} linear combination of ratio $\la$ of the functions $u$ and $v^*$, as defined in \S3.4.
By \cite{L}, $u$ and $v^*$ are $Q^2_-$ functions and hence $\ula\in Q^2_-(\Omega\setminus S_\la)$ and solves (see \cite{CS1})
$$
\begin{cases}
\Deltap \ula \ge 0 \qquad&\text{ in }\oo\setminus\overline\sla\\
\ula =0            \qquad&\text{ on }\bd\oo\\
\ula=1             \qquad&\text{ on }\bd\sla\,.
\end{cases}
$$
Moreover, $|D\ula| \le \La(\oo)$ on $\bd\sla$, thanks to formula (\ref{LS}). Notice that $\Deltap\ula$ can not vanish identically, since $\Omega\setminus\overline\sla$ and $\Omega\setminus\overline\ou$ are not homothetic (see the proof of equality in the Brunn-Minkowski inequality for $\p$-capacity in \cite{CS1}).

Now denote by $w$ the $\p$-capacitary function of $\sla$ with respect to $\oo$ and notice that $\Omega\setminus\overline\sla$  satisfies an uniform interior sphere condition on $\bd\sla$, that is there exists $r>0$ such that for every point $x\in\partial\sla$ there exists a ball $B_r\subset\Omega\setminus\overline\sla$ of radius $r$ with $x\in\bd B_r$ (one can take $r=\mathrm{dist}(\sla,\bd\Omega)/3$ for instance).
Then a careful application of Hopf's Comparison Principle of Tolksdorf [19] 
gives $|Dw|\le \La(\oo)-\ep$ on $\bd\sla$ and hence $(w,\sla)\in \F(\oo,\La(\oo)-\ep)$ which contradicts (\ref{bernconstant}).
\end{proof}


\end{document}